\documentstyle[12pt]{article}
\normalbaselines
\begin{document}
\vbox {\vspace{6mm}}
\def\WP{W(\q,\phi)}
\def\WPa{W(\q,\phi, a)}
\def\FPa{F_a(\phi)}
\def\FP{F(\phi)}
\def\0{\Pi_0}
\def\APR{A^{\phi}(\R)}
\def\OP{\Omega^{\phi}}
\def\AP{A^{\phi}}
\def\hP{\h^{\phi}}
\def\hPR{\hP(\R)}
\def\HPR{\HP(\R)}
\def\HP{H^{\phi}}

\def\ll{\lambda}
\def\o{\omega}
\def\O{\Omega}
\def\W{{\bf W}(\phi)}
\def\r{{\bf r}}
\def\g{{\bf g}}
\def\Wr{{\bf W}(\r,\phi)}
\def\Wb{{\bf W}_{\Phi}}
\def\d{{\bf d}}
\def\R{{\bf R}}
\def\C{{\bf C}}
\def\gC{{\bf g}({\C})}
\def\1C{{\bf g}_1({\C})}
\def\2C{{\bf g}_2({\C})}
\def\dC{{\bf d}({\C})}
\def\dR{{\bf d}(\R)}
\def\gR{{\bf g}(\R)}
\def\lq{\bar\q}
\def\lqq{^{\lq}}
\def\w{{\bf w}}
\def\h{{\bf h}}
\def\t{{\bf t}}
\def\a{\alpha}
\def\E{E_{\a}}
\def\F{E_{-\a}}
\def\H{H_{\a}}
\def\u{{\bf u}}
\def\q{\sigma}
\def\e{\varepsilon}
\def\s{s^{\e}}
\def\th{\theta}
\def\k{{\bf k}}
\def\m{{\bf m}}
\def\b{\beta}
\def\n{{\bf n}}
\def\N{{\bf N}}
\def\p{{\bf p}}
\def\l{{\bf l}}
\def\GC{G(\C)}
\def\GR{G(\R)}
\def\ss{{{\bf g}^*}}
\def\ssC{{\bf g}^*({\C})}
\def\G{\Gamma}
\def\qq{^{\q}}
\def\lh{\overline\h}
\def\lr{\overline\r}
\def\1{\Pi_1}
\def\2{\Pi_2}
\def\AR{A_{\R}}
\begin{center}
{\bf MANIN TRIPLES OF REAL SIMPLE LIE ALGEBRAS. PART 1} 
\end{center}
\begin{center}
A.PANOV\\
Samara State University\\
ul.Akad.Panlova 1,
Samara, 443011,Russia\\
panov@info.ssu.samara.ru
\end{center}
{\bf Abstract.} The article is devoted to the problem of 
classification  of Manin triples
up to weak and gauge equivalence. The case of complex  
simple Lie algebras can be obtained by papers of A.Belavin, V.Drinfel'd,
M.Semenov-Tian-Shanskii. Studing the action of conjugaton on complex
Manin triples, we get the list of real doubles. There exists three types of 
the doubles.
We classify all $ad$-invariant forms on the double compatible with
bialgebra structure.  
We classify the Manin triples $(\gR,W,\gC)$ (case 3 of the doubles) 
up to weak and gauge equivalence.

\begin{center}
{\bf  0.Introduction} 
\end{center}

The problem of classificaton of the Lie-Poisson brackets on Lie groups
leads to the notion of Manin triple .\\ 
{\bf Definition 1}. Let $\g_1,\g_2,\d$ be  Lie algebras over a field $K$
and let $Q$ be a symmetric  nondegenerate bilinear form on $\d$.
A triple $(\g_1, \g_2,\d)$ is called a Manin triple if 
 $Q(x,y)$ is $ad$-invariant and $\d$ is a direct sum of 
maximum isotropic subspaces $\g_1, \g_2$ . 

If $\g = \g_1$ , then $\g_2$ can be identifined with $\ss$. The algebra $\d$ is
 called a double of $\g$. The Lie algebra structure on $\ss$ induces a
Lie coalgebra structure 
on $\g$ and these structures are compatible. In this case we say that
$\g$ is a Lie bialgebra.\\ 
{\bf Definition 2}. We say that two Manin triples 
$(\g,W,\d)$ and $(\g,W',\d)$
are weak equivalent if there exists  an element $a$ in the adjoint group
$D$ of the double $\d$ such that $ W'= Ad_a(W)$.\\
{\bf Definition 3}. We say that two Manin triples 
$(\g,W,\d)$ and $(\g,W',\d)$
are gauge equivalent if there exists  an element $a$ in the adjoint group
$D$ of double  $\d$ such that $ W'=Ad_a(W)$ and $Ad_a (\g) =\g$.\\
Every two gauge equivalent Manin triples are weak equivalent.

One set up the problem of classification  of all Manin triples
up to weak and gauge equivalence in terms of Lie structure of $\g$.
The paper is organized as follows. 
In the section 1 we recall the classification
of Manin triples of simple complex Lie algebras.
In other sections, we study real Manin triples. 
The complexification of a real Manin triple
is a complex Manin triple. We treat the extentions of conjugation $\q$
of the complex Lie algebra $\g$ to the complex double $\d$. We obtain 
the classification of real doubles and the classification 
of the forms $Q$ on the real doubles (see Theorems 2.2, 2.3).  
Further, we study the  complex $\q$-invariant Manin triples.
The main Theorems 4.4 and 4.14 provides the classification of 
Manin triples $(\gR,W,\gC)$ (Case 3 of the double) 
up to weak and gauge equivalence. It turns out that every
class of weak equivalent Manin triples consists of the finite number
of gauge equivalent classes. 
The corresponding $R$-matrices are studied in [CGR], [CH].\\ 
{\bf Acknowlegments}. The author is thankful to 
A.Stolin for the friendly support. His preprint 
[S] was the origin of this paper.\\
Notations:\\
$A_1(\R) = \R[t]/t^2$;
$A_2(\R) = \R[t]/(t^2-1/4)$; $A_3(\R) = \C$;
$A_1(\C) = \C[t]/t^2$;
$A_2(\C) = {\C[t]}/{(t^2-1/4)}.$\\
\begin{center}
{\bf  1. Bialgebra structure of complex simple Lie algebras} 
\end{center}
In this section we recall some facts on bialgebra structure of
simple complex Lie algebras.\\
Let $\g $ be a simple Lie algebra over the field  $\C$ of complex numbers.
Let $G$ be the adjoint group of $\g$. 
Denote by $\Delta $ ( resp. $\Delta ^+$)  
the system of roots (resp.of positive roots) and by 
$\Pi=\{\a_1,\ldots ,\a_n\}$ 
the system of simple roots. Let $K(\cdot,\cdot)$ be a Killing form. 
Consider the Weyl basis : $$H_k=H_{\a_k}, \a_k\in \Pi; E_{\a}, E_{-a}, 
\a\in\Delta^+$$  
$$ [H,\E]=\a(H)\E , \quad [\E ,\F]=H_{\a} ,\quad $$
$$\a(H)=K(H_{\a}, H),\quad K(\E,\F)=1$$
As usial $\h$ is a Cartan subalgebra and $\n^{\pm}$ is up and low
nilpotent subalgebras.\\
{\bf Theorem 1.1}([D],[S]). 
Let $\g$ be a complex simple Lie algebra  and 
let $(\g,W,\d)$ be a complex Manin triple.
 Let $ Q(\cdot, \cdot)$ be nondegenerate, $ad$-invariant form on $\d$
such that $\g$ is an isotropic subspace with respect to $Q$. Then\\
1) $\d$ is isomorphic to $\g\otimes A(\C)$,
where $A(\C)=A_1(\C)$ or $A(\C)= A_2(\C)$;\\
2) $Q(a+bj, a'+b'j) = \lambda(K(a,b')+K(b,a')),$
where  $\lambda \in\C$ and $j$ is an image of t in $A_i(\C)$ for $i=1,2$.\\
Case 1.\\ 
$\d$ is a semidirect sum of $\g$ and commutative subalgebra $\g j$.
Lie algebra $\g$ acts acts on $\g j$ by adjoint action.
The following theorem describes the bialgebra structure. 
{\bf Theorem 1.2}([S]). Let $\g,W,\d$ be a Manin triple and $\d=\g\otimes A_1(\C)$.
Then $$\W=\{a+f(a)j+a^\perp j\vert a\in L, a^\perp \in L^\perp\}$$
Here  $L\subset \g$ is a subalgebra, $L^\perp$ orthogonal
complement with respect to the Killing form, $f:L\mapsto L^*=\g/L^\perp$
is an isomorphism of vector spaces and a 1-cocycle ol L with values is $L^*$.\\
Case 2. \\        
Denote $e=\frac{1}{2}+j$ and $f=\frac{1}{2}-j$.
The elements $e$, $f$ are ortogonal idempotents, 
i.e. $ef=fe=0$ and $e^2=f^2=1$.
In the second case 
$$\d = \g \otimes A_2(\C) = \g e +\g f $$
The simple calculations yield that 
$$Q(ae+bf, a'e+b'f) = K(a,a')-K(b,b')$$ 
Let $\1$, $\2$ be the subsets of $\Pi$ and $\phi:\1\mapsto\2$ is the map
such that:\\
i) $\phi$ is a bijection,\\
ii) $(\phi(\a),\phi(\b))=(\a,\b)$,\\
iii) for every $\a\in\1$ there exists $k$  with the property
$\phi(\a),\ldots,\phi^{k-1}(\a)\in\1$ and $\phi^k(\a)\notin\1$.\\
Notations:\\
$\g_i$ the Lie semisimple subalgebra of $\g$ with the
system of simple roots $\Pi_i$ and Cartan subalgebra $\h_i$.\\
$\phi:\g_1\mapsto\g_2$ isomorphism of semisimple Lie subalgebras defined
via $\phi:\1\mapsto\2$,\\ 
$\p_i^{\pm}$ is upper and lower parabolic subalgebras,
$ \n_i^{\pm}$ is nilpotent radical of $\p_i^{\pm}$,\\
$\r_i=\g_i+\h$ reductive part of $\p_i$', $R_i$ is Lie subgroup of $\r_i$ and 
$R' = R_1\bigcap R_2$\\ 
Let $\lh_1$, $\lh_2$ be the subalgebras of $\h$ such that 
$\lh_i^\perp\subset\lh_i$ and $\h_i\subset\lh_i$. Denote
$\l_i =\lh_i^\perp = Ker K(\cdot,\cdot)\vert_{\lh_i}$, 
$\lr_i = \g_i + \lh_i.$\\
Suppose that  $\phi$ is extended  to an isomorphism
$$\phi :\lr_1/\l_1\mapsto\lr_2/\l_2\eqno (1.1)$$
{\bf Definition 1.3}. The following condition we shall call the condition iv):\\
iv) $\phi$ preserve $K(\cdot,\cdot)$ and has no fixed points in $\h$.\\
( the last means that there is no $x\in\lr_1\bigcap\lr_2$, $x\ne0$ such that 
$\phi(x)=x(mod(\l_2))$.\\
Consider
$$ \W = \l_1 e + \n_1^+ e +\sum_{x\in\lr_1}\C(xe + \phi(x)f)
+ \n_2^ - f +\l_2 f\eqno(1.2)$$
{\bf Theorem 1.4} ([BD],[ST],[S]). 
Let $\phi$ satisfy i), ii), iii), iv). Then \\
1) $(\g, \W, \d)$ is a Manin triple;\\
2) Every Manin triple $(\g, W',\d)$ is gauge 
equivalent to $(\g, \W, \d)$.\\ 
\\
\begin{center}
{\bf  2. Doubles of real simple Lie algebras} 
\end{center}

In this section we are going to study the doubles of real forms of simple 
Lie algebras. The classification of real doubles was obtained in
[CH] by methods of r-matrices.
In our paper we classify the actions 
of conjugation on the complex doubles. This enables us to derive 
the structure of real doubles. We also classify all $ad$-invariant forms 
$Q$ on the double $\dR$ such that $\gR$ 
is an isotropic subspace with respect to $Q$.

It is well known that every real simple Lie algebra is isomorphic to
a real form of simple complex Lie algebra or has the complex structure
[VO], [GG].
The main Theorems of the section are Theorem 2.2 (case of real forms)
and Theorem 2.3 (case of algebras with complex structures).

Let $\g$ be a simple Lie algebra over the field $\C$ and $\gR$
is the real form of $\g$.  
Denote by $\G=\{1,\q : \q^2=1\}$ the Galois group $\C/\R$. There exists
a semilinear action of $\G$ in $\g$ such that $\gR = \g^\G$.

Note that if $(\g_1(\R),\g_2(\R),\dR) $ is a Manin triple over $\R$, then
the  algebras
$\g_1 = \g_1(\R)\otimes\C$,
$\g_2 = \g_2(\R)\otimes\C$,
$\d = \dR\otimes\C$ form a Manin triple over $\C$. 
In the case $\g=\g_1$ we have $\g +\ss=\d$, where $\ss$ is isomorphic
to complex conjugate space $\g(\C)^*$.    
The Galois group acts on each of this algebras and $\gR=\g^\G$,
$\ss(\R)=\ss^\G$, $\dR=\d^\G$.\\
We shall denote by $\lq$ the extention of $\q$ on $\d$. \\
{\bf Lemma 2.1}. Let $\phi$ be a $\C$-linear map of a complex simple Lie 
algebra $\g$ satisfying $[x,\phi(y)] = \phi([x,y])$ for all $x,y\in\g$.
Then $\phi$ is a scalar map (i.e.$\phi(y)=\ll y$ for all $y\in\g$).\\
{\bf Proof}.
 There exists an eigenvector $\phi(v)=\lambda(v)$.
Then $\phi [x,v]= [x, \phi(v)]=[x,\lambda v]=\lambda[x,v]$.
It follows that $[x,v]$ are also the eigenvectors for all $x\in \g$.
Then $[x_1[\ldots [x_k,v]\ldots ]$ are  eigenvectors for all 
$x_1,\ldots ,x_k$ in $\g$. We recall that $\g$ is a simple Lie algebra.
Hence $\phi(y)=\lambda y$. $\Box$\\
{\bf Theorem 2.2}. Let $(\gR, W(\R),\dR$ be a Manin triple for
a real form $\gR$ of a simple complex Lie algebra $\g$. \\
1) There exist the semilinear action of $\G$ in $A(C)$ such, that
$^{\q}(x\otimes a)= ^{\q}x\otimes^{\q}a$ , where 
$ x\otimes a\in \d=\g\otimes A(\C)$;\\
2) The double $\dR$ is isomorphic to  
$\gR\otimes A(\R)$, where $A(\R)$is one of the algebras $A_i(\R)$, $i=1,2,3$;\\
3) Let $ Q_\R(\cdot, \cdot)$ be nondegenerate, $ad$-invariant form on $\dR$
such that $\gR$ is an isotropic subspace with respect to $Q_R$.
Let $K_\R(\cdot,\cdot)$ be a Killing form for $\gR$.
Then
$Q_\R(a+bJ, a'+b'J) = \nu(K_\R(a,b')+K_R(b,a')),$
where  $\nu \in\R$, $J$ is an image of t in $A_i(\R)$ for $i=1,2$
and $J^2=-1$ in $A_3(\R)=\C$.\\
{\bf Proof}. According to Theorem 1.1, the complex double $\d$ is isomorphic to
$\g\otimes A(\C)$, where $A(\C)=A_1(\C)$ or $A(\C)= A_2(\C)$.
Denote by $j$ the image of $t$ in $A_i(\C)$. Then 
$$d=\{x+yj\vert x,y\in\g\}$$
and $j^2 = 0$ or $j^2=\frac{1}{4}$.
Up to definition, ${\lqq x} = {\qq x}$.
There exist the $\C$-linear maps $\phi$ and $\psi$
such that ${\lqq(yj)} = \phi({\qq y}) + \psi({\qq y})j$.
Consider
$${\lqq([x,yj])} ={\lqq([x,y]j)} = \phi({\qq [x,y]}) + \psi({\qq [x,y]})j.$$
On the other hand
$${\lqq([x,yj])} = [{\qq x},{\lqq(yj)]} = 
[{\qq x},\phi({\qq y})] + [{\qq x}, \psi({\qq y})]j$$
Therefore $[x,\phi(y)] = \phi([x,y])$ and $[x,\psi(y)] = \psi([x,y])$.
According to Lemma 2.1, $\phi(y)=\ll y$ and $\psi(y)=\mu y$. 
Denote by $J$ the element $\ll+\mu j$ in $A(\C)$. We have
$$ {\lqq(yj)}= {\qq(y)}J\eqno (2.1)$$
For every two elements $y, y'\in\g$ there hold
$${\lqq([yj, y'j])} = {\lqq([y, y']j^2)} = j^2{\qq([y,y'])} = 
j^2[{\qq y}, {\qq y'}]$$
On the other hand,
 $$ {\lqq([yj, y'j])} = [{\lqq (yj)}, {\lqq (y'j)}] = 
[{\qq y}J, {\qq y'}J] = J^2[{\qq y}, {\qq y'}]$$
We have 
$$ J^2 = j^2 \eqno (2.2)$$
The following three case take place.\\
Case 1. $A(\C) = A_1(\C)$. Since $j^2=0$ then $J^2=0$. Hence, $\ll=0$ and
$J = \m j$. Substitute this in (2.1)
$$ {\lqq (yj)} = \mu({\qq y})j$$
The $\lq^2 = 1$ implies $ \overline{\mu} \mu =1$.
There exists
$\nu\in\C$ such that $\mu=\overline{\nu}^{-1}\nu$. Denote $J_1=\nu j$.
We have  
$${\lqq(yJ_1)} = {\lqq(y\nu j)} = \overline{\nu}({\qq y})\mu j = {\qq y}\nu j
= {\qq y} J_1$$
Every element in $\d$ have the form $x+yJ_1$, $J_1^2=0$.
The element is $\lq$-invariant if and only if ${\qq x} = x$ and ${\qq y} = y$.
Then $\dR = \d^{\Gamma} = \gR+\gR J_1 = \gR\otimes A_1(\R)$.
Case 2. $j^2=\frac{1}{4}$, $J=j$. Then ${\lqq (yj)}= {\qq y}j$
and ${\lqq (x+yj)} = {\qq x} + {\qq y}j$.
Therefore, $\dR = \d^{\Gamma} = \gR+\gR j = \gR\otimes A_2(\R)$.\\
Case 3. $j^2=\frac{1}{4}$, $J= -j$. Then ${\lqq (yj)}= -{\qq y}j$
and ${\lqq (x+yj)} = {\qq x} - {\qq y}j$.
Therefore,  $\dR = \d^{\Gamma} = \gR+\gR ij = \g$. 
This proves the statements 1) and 2) of the Theorem.

Let $ Q_\R(\cdot, \cdot)$ be nondegenerate, $ad$-invariant form on $\dR$
such that $\gR$ is an isotropic subspace with respect to $Q_\R$.
We extend $Q_\R$ on $\d$ and get the form $Q$, satisfying
the conditions of the 2) from Theorem 1.1.
 For the elements $a,b,a',b'\in \gR$ we have
$$Q_\R(a+bj, a'+b'j) = \lambda(K(a,b')+K(b,a'))$$
Note that $K_\R(x,y)=K(x,y)$ for $x,y\in\gR$. 
This proves statement 3) in Cases 1 and 2.
In the case 3 put $J=2ij$ and
$Q_\R(a+bJ, a'+b'J) = 2i\lambda(K_\R(a,b')+K_\R(b,a'))$
$\Box$.

Our next goal to consider  real simple Lie algebras with the 
complex structure.\\
{\bf Theorem 2.3}. Let $\g$ be a complex simple Lie algebra
and $\g_\R$ is realification of $\g$. Let $(\g_\R, W_\R, \d_\R)$
be a Manin triple. We assert that \\
1) The double $\d_\R$ is isomorphic to $\d_\R\otimes A_i(\R)$ for $i=1,2$.\\
2) Let $ Q_\R(\cdot, \cdot)$ be nondegenerate, $ad$-invariant form on $\d_\R$
such that $\g_\R$ is an isotropic subspace with respect to $Q_R$.
Let $K_\R(\cdot,\cdot)$ be a Killing form for $\g_\R$.
Then
$Q_\R(a+bj, a'+b'j) = \nu(K_\R(a,b')+K_\R(b,a')),$
where  $\nu \in\R$, $j$ is an image of t in $A_i(\R)$ for $i=1,2$.\\
{\bf Proof}. If $(\g_\R, W_\R, \d_\R)$ is a Manin triple then
$(\g_\R\otimes \C, W_\R\otimes\C, \d_\R\otimes \C)$ is also a Manin triple.
Let $\C=\{a+bk\vert k^2=-1; a,b\in\R\}$. We shall denote the complex structure
on $\g_\R$ by $i$. Then the elements
$p=\frac{1+ik}{2}$ and $q=\frac{1-ik}{2}$ form the system of orthogonal 
idempotents in $\C\otimes \C$ and $\g_\R \otimes\C = \g p\oplus\g q$. Here
$\g p$ and $\g q$ are complex algebras isomorphic to $\g$.
The conjugation $\q$ acts on $\g_\R \otimes\C $ by ${\qq (xp)} = xq$ and
${\qq (xq)} = xp$. The $\g_\R = \{xp + xq\}$.

The double $\d_\R\otimes\C$ of $\g_\R\otimes\C$ is isomorphic to the 
direct some of the doubles
$\d_\R\otimes\C = \d_1 \oplus\d_2$ where $\d_1$ is a double
of $\g p$ and $\d_2$ is a double of $\g q$.
Accoding to the Theorem 1.4 , $\d_1 = \g p + \g pj_1$,
$\d_2 = \g_q + \g qj_2$ where $j_1^2$ and $j_2^2$ are separately equal to
0 or $\frac{1}{4}$.
Denote by $\lq$ the conjugation on $\d_\R\otimes\C$ over $\g_\R$.
Recall that the retriction of $\lq$ on $\g_\R\otimes\C$ coinsides with $\q$.
By definition,
$${\lqq(ypj_1)} = \th_1(y)p + \th_2(y)q + \eta_1(y)pj_1 + \eta_2(y)qj_2$$
where $\th_1, \eta_1 $ are $\C$-antilinear  and 
$\th_2, \eta_2 $ are $\C$-linear maps $\q\mapsto\g$.
We have
$${\lqq([xp, ypj_1])} = {\lqq([x, y]pj_1])} = 
\th_1([x,y])p + \th_2([x,y])q + \eta_1([x,y])pj_1 + \eta_2([x,y])qj_2$$
On the other hand,
${\lqq([xp, ypj_1])} = [{\lqq(xp)}, {\lqq(ypj_1)}] = 
[xq , \th_1(y)p + \th_2(y)q + \eta_1(y)pj_1 + \eta_2(y)qj_2] = 
[x, \th_2(y)]q + [x, \eta_2(y)]qj_2$.
Comparing this two formulas, we see $\th_1([x,y]) = \eta_1([x,y]) =0$
(therefore, $\th_1=0$ and $\eta_1=0$) 
and
$[x,\th_2(y)] = \th_2([x,y])$, $[x,\eta_2(y)] = \eta_2([x,y])$.
By Lemma 2.1,  $\th_2(y) = \ll y$ and $\eta_2(y) = \mu y$.
Thus, ${\lqq(ypj_1)} = \ll yq + \mu yqj_2$.
Similarly, ${\lqq(yqj_2)} =\ll ' yp + \mu ' ypj_1$.\\
Denote
$$\left \{ \begin{array}{rcl}J_1&=& \ll'+\mu'j_1\\
J_2& = &\ll+\mu j_2\end{array}\right.\eqno (2.3)$$
We have
$$\left \{ \begin{array}{rcl}{\lqq(ypj_1)} &=&yqJ_2 \\
{\lqq(yqj_2)} & = &ypJ_1\end{array}\right.\eqno (2.4)$$
Compare the left and right sides of the following equalities:
$${\lqq([ypj_1, y'pj_1])} = {\lqq([y, y']pj_1^2)} = [y, y']qj_1^2$$
$$ {\lqq([ypj_1, y'pj_1])} = [ {\lqq(ypj_1)}, {\lqq(y'pj_1)}] =
[yqJ_2, y'qJ-2] = [y, y']qJ_2^2$$
It follows $$j_1^2 = J_2^2$$ 
 Similarly, $$j_2^2 = J_1^2$$
Consider the following two cases i) and ii).\\
Case i): $j_1^2=0$. Then $J_2=0$. If $j_2^2=\frac{1}{4}$ then the algebra
$\C +\C j_2$ has no nilpotents. A contradiction.

Then $j_2^2 = 0$ and $J_2 = \mu j_2$. 
We have 
$$\left \{ \begin{array}{rcl}{\lqq(ypj_1)} &=&\mu yqj_2 \\
{\lqq(yqj_2)} & = &\mu'ypj_1\end{array}\right.$$
The $\lq^2=1$ yields $\overline{\mu}\mu' = 1$.

Denote $I_1 = \a j_1$, $I_2 = \b j_2$ for $\a,\b\in\C$.
Then  
$$\left \{ \begin{array}{rcl}{\lqq(ypI_1)} &=&\mu\bar\a\b^{-1} yqI_2 \\
{\lqq(yqI_2)} & = &\overline\mu^{-1}\bar\b\a^{-1}ypI_1\end{array}\right.
\eqno (2.5)$$
 Choose $\a,\b$ such that $\mu\bar\a\b^{-1} = 1$ then
$\overline\mu^{-1}\bar\b\a^{-1} = 1$ and
$$\left \{ \begin{array}{rcl}{\lqq(ypI_1)} &=& yqI_2 \\
{\lqq(yqI_2)} & = &ypI_1\end{array}\right.\eqno (2.6)$$
The element $xp+xq+ypI_1+y'qI_2$ is $\lq$-invariant
if $x=x'$ and $y=y'$.
The double $\d_\R $ is isomorphic to 
$\g_\R\otimes A_1(\R) = \g_\R + \g_R I $, $I^2=0$.\\
Case ii):  $j_1^2 = \frac{1}{4}$ . Then
$$\frac{1}{4} = J_2^2 = (\ll+\mu j_2)^2 = \ll^2+2\ll\mu j_2+\mu^2j_2^2$$
This leads to
$$\left \{ \begin{array}{rcl}\ll^2+\frac{1}{4}\mu^2 &=&\frac{1}{4} \\
\ll\mu & = &0\end{array}\right.\eqno (2.7)$$
If $\mu=0$, then 
${\lqq(ypj_1)}\in\g q = {\lqq(\g p)}$. A contradiction.
Hence $\ll=0$ and $\mu=\pm 1$.
Substituting this to (2.4), we have
$$\left \{ \begin{array}{rcl}{\lqq(ypj_1)} &=&\pm yqj_2 \\
{\lqq(yqj_2)} & = &\pm ypj_1\end{array}\right.\eqno (2.8)$$
Here $+$ and $-$ appears simultaneously. In the case of $-$
we change the generator $j_2$ to $-j_2$. Finally,
$$\left \{ \begin{array}{rcl}{\lqq(ypj_1)} &=& yqj_2 \\
{\lqq(yqj_2)} & = &ypj_1\end{array}\right.\eqno (2.9)$$
The element $xp+xq+ypI_1+y'qI_2$ is $\lq$-invariant
if $x=x'$ and $y=y'$.
The double $\d_\R $ is isomorphic to 
$\g_\R\otimes A_2(\R) = \g_\R + \g_\R I $, $I^2=\frac{1}{2}$.
This proves the first statement of the Theorem.

Let $ Q_\R(\cdot, \cdot)$ be nondegenerate $ad$-invariant form on $\d_\R$
such that $\g_\R$ is an isotropic subspace with respect to $Q_R$.
We extend $Q_\R$ to nondegenerate $ad$-invariant $\C$-form  $Q$ on 
 $\d_\R = \d_1+\d_2$. The algebra $\g_\R = \g p\oplus\g q$ is
isotropic subspace with respect to $Q$.
Then $Q(\d_1,\d_2)=0$ and  the restriction of $Q$ on $\d_1$, $\d_2$
can be calculated by the statement of Theorem 1.1 with 
constants $\ll_1$, $\ll_2$.

 For all 
$a,b,a',b'\in\g_\R$ we obtain
$Q_\R(a+bj, a'+b'j) = Q(ap+aq+bpj+bqj, a'p+a'q+b'pj+b'qj) = 
\ll_1 (K(ap,b'p)+K(bp,ap')) +\ll_2 (K(aq,b'q)+K(bq,aq'))$.
The map $a\mapsto aq$ (resp.$a\mapsto ap$) is $\C$-linear
(resp. $\C$-antilinear) isomorphism of $\g$ onto $\g q$ (resp.$\g p$).
Therefore, $ K(ap,b'p) = \overline{K(aq, b'q)} = \overline{K(a,b')}$.
Similarly, $ K(bp,a'p) = \overline{K(bq, a'q)} = \overline{K(b,a')}$.
Since, $Q_\R(a+bj, a'+b'j) \in\R$ then $ \ll_1=\ll_2=\nu\in\R$.
We have 
$Q_\R(a+bj, a'+b'j) = \nu (K(a,b')+\overline {K(a,b')} +
K(b,a')) + \overline{K(b,a')}) = \nu( K_\R(a,b') + K_\R(b,a')) $.
$\Box$\\
\begin{center}
{\bf  3. On real forms of complex simple Lie algebras} 
\end{center}

As above $\g = \gC$ is a simple Lie algebra over the field  $\C$
and $G=G(\C)$ is the adjoint Lie group of $\g$.

We shall recall some statements on real forms.
Every compact form of $\g$ is $Ad_{G(\C)}$-conjugates
to following standart compact form :
$$\u = \sum_{\a\in\Pi}\R iH_{\a} + \sum_{\a\in\Delta^+} \R(\E -\F)
+ \sum_{\a\in\Delta^+} \R i(\E +\F)$$
Denote by $\tau$ the conjugation of $\g$ over $\u$ :
$$ \tau (\E) = -\F, \quad\tau(\F) = -\E,\quad\tau(\H) = -\H$$
The algebra 
$$\t=\sum_{\a\in\Pi}\R i\H$$ is a maximum commutative subalgebra in $\u$.

We consider a real form $\gR$ of $\g$ with conjugation $\q$.
Denote $$G(\R) = G(\C)^\Gamma= \{g\in\G(\C)\vert {\qq g}= g\}$$ 
There exists a compact form with cojugation commutating with 
$\q$ [VO,GG]. We assume that this compact form coinsides with $\u$,
i.e. $\q\tau = \tau\q$. Denote $\theta = \q\tau$. The correspodence
$\g(\R)\quad\mapsto\quad\theta$ is 1-1 between the set of real 
forms of $\g$ and  automorphisms of the second order of $\g$.
We can choose $\g(\R)$ in the  $Ad_{G(\C)}$-class such that
$$ \theta = \s exp( ad_{h'}),\eqno (3.1)$$ 
where $s$ is the automorphism of the second order of Dynkin's diagram
( for $A_n, D_n, E_6$), $\e = 0,1$ and $h'\in\h$  commutates 
with $s$.
The $\g(\R)$ is called an inner form if  $\e=0$ and  is called
an outer form if $\e=1$.\\
$$\k = \{x\in \u\vert \quad\th(x)=x\}$$
$$\m = \{x\in\u\vert \quad\th(x)=-x\}$$
Then $\u = \k + \m$ and $\g(\R) = \k + i\m$.\\
Our aim is to classify all Manin triples for $\g(\R)$.
We see in the section 2 that there are three possibilities 
for the double $\dR$ of $\gR$.\\
\\
\begin{center}
{\bf 4. On Manin triples $(\gR,W,\gC)$}
\end{center}
In this section we shall consider the Manin triples $(\gR,W,\dR)$ with
$\dR=\gC$ (i.e. Case 3).
Then $(\gR\otimes\C, W\otimes\C, \dR\otimes\C)$ is Manin triple over $\C$.
Recall that $\gR\otimes\C = \g$ and $\d = \dR\otimes\C = \g e + \g f$.
The conjugation $\q$ acts on $\g$ and $\d$.
In this case $\q e = f$, $\q f = e$.\\
According to Theorem 1.4,  $W\otimes\C$ coinsides with some $\q$-invariant
subalgebra $Ad_g\W$, where $g\in G$.
We are going to classify all $\q$-invariant  subalgebras in the orbit
$Ad_g\W$ and to get there real forms $ W=\W^\Gamma$.
Consider $$ ^\q(Ad_g\W) =Ad_g\W\eqno (4.1)$$
It follows
$${^\q(\W)} = Ad_p\W\eqno (4.2)$$
where $p={^\q g^{-1}} g$. Note that $^\q p = p^{-1}$.  
Replacing $\W$ by (1.2) in (4.1), we get the equality of two  subalgebras
$$\begin{array}{c}^\q\l_1f + ^\q \n_1^+ f +
\sum_{x\in\lr_1}\C(^\q xf + ^\q(\phi(x))e)
+^\q \n_2^ - e + {^\q\l_2}e =   
Ad_p(\l_1)e +\\ 
Ad_p(\n_1^+)e+
\sum_{x\in\lr_1}\C(Ad_p(x)e + Ad_p(\phi(x))f)
+ Ad_p(\n_2^ -) f + Ad_p(\l_2)f\end {array}\quad  \eqno (4.3)$$
{\bf Lemma 4.1.} Let $\p$ , $\p$ be paraboic subalgebras in $\g$
and $\n^+\subset \p,\p'$.
Suppose that $\p' = Ad_g \p$ for some $g\in G$. Then $\p=\p'$.\\
{\bf Proof.}   The $ Ad_g\n^+$ and $\n^+$ are maximum nilpotent subalgebras in 
$\p'$. Let $P'$ be the adjoint group of $\p'$.
There exists $\g_0\in P'$ such that $Ad_{g_0} Ad_g\n^+ = \n^+$.
Then $Ad_{g_og}\n^+ = \n^+$ and $g_0g$ belongs to the Borel
subgroup $B^+$. Denote $g_og = b$. Then
$$\p' = Ad_{g_o}\p' = Ad_{g_o}Ad_g\p = Ad_{g_og}\p = Ad_b\p = \p$$ 
$\Box$\\
{\bf Lemma 4.2}. We save the notation of section 1. Let (4.3) holds.
Denote $\m_1=\g_1\bigcap\n_2^-$, $\m_2=\g_2\bigcap\n_1^+$.
We assert that \\
1) The conjugation $\q$ takes $\l_1$, $\r_1$, $\lr_1$, $\g_1$, $\m_1$ to
$\l_2$, $\r_2$, $\lr_2$, $\g_2$, $\m_2$
and $\l_2$, $\r_2$, $\lr_2$, $\g_2$, $\m_2$ to 
$\l_1$, $\r_1$, $\lr_1$, $\g_1$, $\m_1$.
This defines the maps 
$\q :\lr_1/\l_1\mapsto\lr_2/\l_2$,
$\q :\lr_2/\l_2\mapsto\lr_1/\l_1$\\
2) $p\in P_1^+\bigcap P_2^-$ and $p= {^\q m^{-1}qm}$ where
$m\in\exp(\m_2)$, ${^\q m} \in exp(\m_1)$, ${^\q q} = q^{-1}$, $q\in R'$; \\
3) Denote $ p_1 = {^\q m^{-1}q}$ and $p_2 = qm$. Then 
$Ad_{p_1}:\lr_1/\l_1\mapsto\lr_1/\l_1$,
$Ad_{p_2}:\lr_2/\l_2\mapsto\lr_2/\l_2$
The map $^\q Ad_{p_2}\phi$ takes $\lr_1/\l_1$ to $\lr_1/\l_1$ and
$$(^\q Ad_{p_2}\phi)^2(x)= x\eqno (4.4)$$
for $x\in\lr_1/\l_1$.\\
{\bf Proof}. 
Consider the intesection of left and right sets with $\g e$ and $\g f$.\\
We have 
$$\left \{ \begin{array}{rcl}^\q(\l_1+\n_1^+)&=& Ad_p(\l_2+\n_2^-)\\
^\q(\l_1+\n_2^-)& = &Ad_p(\l_1+\n_1^+)\end{array}\right.\eqno (4.5)$$
Compairing the radicals of subalgebras, we have

$$\left \{ \begin{array}{rcl}^\q\n_1^+&=& Ad_p\n_2^-\\
^\q\n_2^-& = &Ad_p\n_1^+\end{array}\right.\eqno (4.6)$$

$$\left \{ \begin{array}{rcl}^\th\n_1^-&=& Ad_p\n_2^-\\
^\th\n_2^+& = &Ad_p\n_1^+\end{array}\right.\eqno (4.7)$$

$$\left \{ \begin{array}{rcl}\s (\p_1^-)&=& Ad_p\p_2^-\\
\s(\p_2^+)& = &Ad_p\p_1^+\end{array}\right.\eqno (4.8)$$

Applying  Lemma 4.1 to pairs of parabolic subalgebras
$\s(\p_1^-), \p_2^-$ and $\s(\p_2^+), \p_1^+$, we get

$$\left \{ \begin{array}{rcl}Ad_p\p_2^-&=& \p_2^-\\
Ad_p\p_1^+& = &\p_1^+\end{array}\right.\eqno (4.9)$$

$$\left \{ \begin{array}{rcl}\s(\p_1^-)&=& \p_2^-\\
\s(\p_2^+)& = &\p_1^+\end{array}\right.\eqno (4.10)$$

The (4.9) implies $p\in P_1^+\bigcap P_2^-$. It follows form
(4.10) that
$$\left \{ \begin{array}{rcl}^\q(\p_1^+)&=& \p_2^-\\
^\q(\p_2^-)& = &\p_1^+\end{array}\right.\eqno (4.11)$$

$$\left \{ \begin{array}{rcl}^\q\n_1^+&=& \n_2^-\\
^\q\n_2^-& = &\n_1^+\end{array}\right.\eqno (4.12)$$

$$\left \{ \begin{array}{rcl}\s(\r_1)&=& \r_2\\
\s(\r_2)& = &\r_1\end{array}\right.\eqno (4.13)$$

The equalities (4.13) yield

$$\left \{ \begin{array}{rcl}^\q(\g_1)&=& \g_2\\
^\q(\g_2)& = &\g_1\end{array}\right.\eqno (4.14)$$

$$\left \{ \begin{array}{rcl}^\q(\m_1)&=& \m_2\\
^\q(\m_2)& = &\m_1\end{array}\right.\eqno (4.15)$$

The element $p$ can be uniquely decomposed in the form
$$p = m_1qm_2\eqno (4.16),$$
where $m_1\in exp(\m_1)$, $m_2\in exp(\m_2)$, $q\in R'$.\\
Denote $p_1=m_1q\in R_1$, $p_2= gm_2\in R_2$. Then $p = p_1m_2 = m_1p_2$.
Since  $^\q p = {^\q m_1}{^\q q}{^\q m_2}$ and 
$p^{-1} = m_2^{-1}q^{-1}m_1^{-1}$,
then $ {^\q m_1} = m_2^{-1}$ and ${^\q q} = q^{-1}$. We have
$$ p = {^\q m^{-1}}qm\eqno (4.17)$$
where $m = m_2\in\m_2$. We see $p_1 = {^\q m^{-1}}q$ and 
$$p_1^{-1} = q^{-1}{\q m} = {^\q({^\q q^{-1}m})} = {^\q(qm)} = {^\q p_2}$$
We have 
$$\left \{ \begin{array}{rcl}^\q p_1&=& p_2^{-1}\\
^\q p_2& = &p_1^{-1}\end{array}\right.\eqno (4.18)$$
Since $\g_2\bigcap\l_2 = 0$, then $Ad_{p_2}\l_2 = \l_2$.
Therefore,
$$Ad_p(\l_2+\n_2^-) = Ad_p\l_2 + Ad_p\n_2^- = 
Ad_{m_1}\l_2 + \n_2^- =\l_2 +\n_2^-$$ 
Substituting this to (4.5) and to (4.3), we have
$$\left \{ \begin{array}{rcl}^\q(\l_1+\n_1^+)&=& \l_2+\n_2^-\\
^\q(\l_1+\n_2^-)& = &\l_1+\n_1^+\end{array}\right.\eqno (4.19)$$
$$\begin{array}{c}^\q\l_1f + ^\q \n_1^+ f +
\sum_{x\in\lr_1}\C(^\q xf + ^\q(\phi(x))e)
+^\q \n_2^ - e + {^\q\l_2}e =   
\l_1e +\\ 
\n_1^+e +
\sum_{x\in\lr_1}\C(Ad_p(x)e + Ad_p(\phi(x))f)
+ \n_2^- f + \l_2f\end {array}\quad  \eqno (4.20)$$
We stress that\\
 $Ad_px= Ad_{p_1}x(mod(\n_1^+))$, 
$Ad_{p_1}x\in\lr_1\subset\r_1$
and $Ad_p\phi(x) =Ad_{p_2}\phi(x)(mod(\n_2^+))$, $Ad_{p_2}\phi(x) \in\lr_2
\subset\r_2$.
It follows from (4.13) and  $^\q\h=\h$ that 
$$\left \{ \begin{array}{rcl}^\q(\l_1)&=& \l_2\\
^\q(\l_2)& = &\l_1\end{array}\right.\eqno (4.21)$$
The (4.20) implies
$$ 
\sum_{x\in\lr_1}\C(^\q xf + ^\q(\phi(x))e)
\subset \l_1e + \n_1^+e +
\sum_{x\in\lr_1}\C(Ad_{p_1}(x)e + Ad_{p_2}(\phi(x))f)
+ \n_2^ - f + \l_2f \eqno (4.22)$$
Since $p_1=h_1g_1$ with $h_1\in H=exp(\h)$, $g_1\in G_1=exp(\g_1)$,
then 
$$\left \{ \begin{array}{rcl}Ad_{p_1}\lr_1&=& \lr_1\\
Ad_{p_2}\lr_2& = &\lr_2\end{array}\right.\eqno (4.23)$$
$$\left \{ \begin{array}{rcl}Ad_{p_1}\l_1&=& \l_1\\
Ad_{p_2}\l_2& = &\l_2\end{array}\right.\eqno (4.24)$$
The (4.22) implies that for every $x\in\lr_1$ 
there exists $y\in\lr_1$ such that

$$\left \{ \begin{array}{rcl}^\q(x) &=&Ad_{p_2}\phi(y)(mod(\l_2))\\
^\q(\phi(x))& = &Ad_{p_1}y(mod(\l_1))\end{array}\right.\eqno (4.25)$$
It follows
$$\left \{ \begin{array}{rcl}^\q(\lr_1)&=& \lr_2\\
^\q(\lr_2)& = &\lr_1\end{array}\right.\eqno (4.26)$$
Finally we have
$$Ad_{p_1}^{-1}{^\q(\phi(x))} = \phi^{-1} Ad_{p_2}^{-1}{^\q(x)}\eqno (4.27)$$
for $x\in\lr_1/\l_1$.\\
Recall that  $p_1^{-1}= {^\q p_2}$.
The (4.27) can be rewritten in the equivalent form
$(^\q Ad_{p_2}\phi)^2(x)= x$. 
for $x\in\lr_1/\l_1$. This proves (4.4). $\Box$\\
{\bf Definition 4.3}. The following condition we shall call the condition v) for
$\phi :\lr_1/\l_1\mapsto\lr_2/\l_2$:\\
v) $\qq\lr_1=\lr_2$, $\qq \l_1=\l_2$ and 
 $${\qq(\phi(x))} = \phi^{-1}({\qq x})\eqno(4.28)$$ for $x\in\lr_1/\l_1.$\\
Note that:\\
1) In the case $\1=\2=\emptyset$ the condition v) means 
$\qq(\phi(x)) = \phi^{-1}(\qq x)$ for $x\in\lh_1$;\\
2) the condition v) can be rewritten in the form
$({^\q\phi})^2(x)= x$ for $x\in\lr_1/\l_1$.\\
3) If we put $x=H_{\a}$, then ${\qq(\phi(H_{\a}))} = \phi^{-1}({\qq H_{\a}})$ 
yields ${^s\phi(\a)}=\phi^{-1}s(a)$ for all $\a\in\1$.\\
Notations:\\
$ \h(\R) = \gR\bigcap\h$;\\
$H(\R) = exp (\h(\R))$;\\
$\AR = \{a\in H(\R)\vert\quad a^2 = 1\}$;\\
$\AR' = \{a\in\AR\vert\quad\exists g\in G, {\qq g^{-1}}g=a\};$\\
$\O(\AR) = \{g\in G\vert\quad{\qq g^{-1}}g \in A\}$;\\
Note that $\h(\R) = \t$ for inner real forms.\\
The group $G(\R)$ acts on $\O(\AR)$  by $g\mapsto rg$, $r\in G(\R)$.
Denote by $\O$ the set of representatives of $G(\R)$- orbites in $\O(\AR)$.
The set $\O$ is finite.\\
{\bf Theorem 4.4} Let $\gR$ be an inner form of $\g=\gC$ . The conjugation
of $\gR$ as in (3.1) with $\e=0$. 
The map $\phi $ as in Theorem 1.4. Then\\
1) $\W $ is $\q$-invariant if and only if  $\1=\2=\emptyset$ and
$\phi$ satisfies v).\\
2) There exists a $\q$-invariant subalgebra in $Ad_{G(\C)}\W$ iff
$\1=\2=\emptyset$ and $\phi$ satisfies v). 
Every $\q$-invariant subalgebra in 
$Ad_{G(\C)}\W$ has a form $Ad_r Ad_\o\W$, where $r\in G(\R)$, $\o\in\O$.
If $\gR=\u$ is a compact form of $\g$ then every $\q$-invariant subalgebra in 
$Ad_{G(\C)}\W$ has a form $Ad_r W$, where $r\in\GR$.\\
{\bf Proof.}
Statement 1) can be proved by direct calculations. We shall prove 2).

Let $Ad_g\W$ be $\q$-invariant subalgebra.
It follows from (4.13) that $\r_1=\r_2$ for inner forms.
Consequently, $\g_1 = \g_2$ and $\1 = \2$.
 This contrudicts to the condition iii) in the definition of $\phi$. 
Hence, $\1 = \2 =\emptyset$, $\r_1=\r_2=\h$ and $\lr_i = \lh_i$.
Using Lemma 4.2, we get $p= p_2\in P_1^+ \bigcap P_2^-= B^+\bigcap B^- = H$,
where $H=exp (\h)$.\\
The equality (4.4) take the form (4.28): 
$^\q(\phi(x)) = \phi^{-1} (^\q x)$
for $x\in\lh_1/\l_1$.
 
Since $p\in H$, then $p=exp(v)$, where
$v = \sum_{k=1}^n \nu_k H_k.$
The linear operator $Ad_p$ is diagonalizable. 
Denote $p_\a=exp(a(v))$, $\a\in\Pi$.
We shall idenify $p$ with system of numbers $(p_\a)_{\a\in\Pi}$.

Note the condition $^\q p=p^{-1}$ is equvalent to $\overline{p_\a}=p_\a$,i.e.
$p_\a\in\R$ and $p = {^\q h^{-1}}h$, $h\in H$ is equvalent to $p_\a>0$.\\
Note that the subalgebra 
$\AR$ coinsides with $\{a\in H\vert  a_\a=\pm 1 \}.$
Then  there exist $h\in H$ and $a\in\AR$ such that $p$ = ${^\q h^{-1}a h} $ .
Hence, ${^\q g^{-1}g}$ = ${^\q h^{-1}a h}$. We see that $a\in \AR'$
We denote by $\o$ in $\O$ such that
$^\q \o^{-1} \o = a$. 
 $$ ^\q (gh^{-1})^{-1} gh^{-1} = {^\q \o^{-1}}\o$$
 $$ ^\q (gh^{-1}\o^{-1})^{-1} gh^{-1}\o^{-1} = 1$$
Therefore, $gh^{-1}\o^{-1} = r \in G(\R)$
and $g=r\o h$.
We have
$$Ad_g\W = Ad_rAd_{\o}Ad_h\W =  Ad_r Ad_{\o}\W $$
This proves 2). 

Now we shall study the case of compact real form $\gR=\u$.\\
Consider the scalar product $<x,y>_{\tau} = -K(\tau x,y)$.
Then $^\q (Ad_g)^{-1} $ coinsides with the conjugating linear operator 
$(Ad_g)^*$.
Then $Ad_p =Ad _{^\q g^{-1}g} =(Ad_g)^*Ad_g$ is positive definite
operator.
On the other hand, the operator $Ad_p, p\in H $  is positive definite
if and only if $exp(\a(v)) > 0$
for all $\a\in \Delta$.\\
It follows $$\a(v)= \sum_k^n\nu_k\a(H_k)\in\R$$
Then $\nu_k\in\R$ for all $k$ and, consequently, $p={^\q h^{-1}}h$. 
That is $\o =1$.
Finally, we get $g=rh$, $r\in G(\R)$ and
$$Ad_g\W = Ad_r Ad_h\W =  Ad_r \W $$
$\Box$

The above Theorem 4.3 gives  the classification of Manin triples 
in the following terms. 
If the subalgebra $Ad_g\W$ is $\q$-invariant, then consider
the $\R$-algebra of invariants
$W(\R) = (Ag_g\W)^\G$ in $\dR=(\dC)^\G$. Recall that 
$\dR$ is isomorphic to $\gC$ (Case 3 of the double structures).  
Using the isomorphism of $\dR$ to $\gC$,
we are going to give more precise  description for $W(\R)$.  

The double algebra $\dR = (\dC)^{\Gamma} $ consists of the
elements $\tilde x = xe+{^\q x}f.$ 
The map $x\mapsto\tilde x$  is an $\R$-isomorphism of $\gC$ to $\dR$.
Recall that $$Q(ae+bf,a'e+b'f)=K(a,a')-K(b,b')$$ is a nondegenerate, $ad$-invariant
bilinear form on $\dC$.
Restricting $Q$ on $\dR$, we get 
$$Q_{\R}(x, y) = K(x,y)-K(^\q x,^\q y) = K(x,y) - \overline{K(x,y)} = 
2iIm(K(x,y))$$ i.e the nondegenerate, 
$ad$-invariant bilinear form on $\dR$. 
Lie algebra $\gR$ is an isotropic subspace with 
respect to $Q_{\R}$.\\
{\bf Theorem 4.5}. Let $\gR$ be the inner form of $\gC$.
Let $\Phi$ be a $\R$-subspace in $\h$ such that $\Phi$ is an  isotropic with
respect to  $Q_\R(\cdot, \cdot)$ and $\Phi\oplus\t =\h$.
Denote $\Wb$ the sudalgebra $\Wb = \Phi \oplus \n^+$. Then:\\
1) $(\gR,\Wb,\gC)$ is a Manin triple;\\
2) Every Manin triple $(\gR,W,\gC)$ is weak equivalent to
 some triple $(\gR, \Wb, \gC)$;\\
3) Every Manin triple $(\gR,W,\gC)$ is gauge equivalent to the Manin triple
$$(\gR,Ad_\o\Wb,\gC),$$ where $\o\in\O$. If $\gR=\u$ is a compact real form, 
then every Manin triple $(\gR,W,\gC)$ is gauge equivalent to the Manin triple
$(\gR,\Wb,\gC)$.\\
{\bf Proof.} 
The proof of the Theorem is consists of the several steps.
On the 1)-3) steps we define the map $\phi$ and prove that it satisfy i)-v)
(Section 1 and Definition 4.3) if and only if $\Phi$ satisfy the assumptions
of the Theorem. On the 4) step we show that $\W$, constructed via $\phi$
(Theorem 1.4), is $\q$-invariant and $\W^\G = \Phi\oplus\n^+$.
On the last 5) step we see that $(Ad_g\W)^\G = Ad_g(\Phi\oplus\n^+)$
for $g=r\o$. The application of Theorem 4.3 proves the Theorem.\\
1) For every $\R$-subspace $V$ in $\h$ we denote $\C V$ the complex subspace 
generated by $V$. Denote by $\l_1$ the greatest complex subspace in the
real subspace $\Phi$, $\lh_1=\C \Phi$, $\lh_2 = \C{^\q\Phi}$, 
$\l_2 = {^\q \l_2}$. 
Note that $\Phi/\l_1$ is $\R$-form of $\lh_1/\l_1$.
There exists the conjugation $\ll$ on $\lh_1/\l_1$ such that
$$\Phi/\l_1= \{x\in\lh_1/\l_1\vert {^\ll x} = x\}$$
The conjugation $\q$
maps $\lh_1/\l_1$ to $\lh_2/\l_2$. By $\phi = \q\ll$ we denote
$\C$-linear map $\lh_1/\l_1$ to $\lh_2/\l_2$.
Since $\q^2=1$, then $\q\phi = \ll$ and $(\q\phi)^2=1$. The map $\phi$
satify v).\\
The conditions i)-iii) are true for $\phi$ with $\1=\2=\emptyset$.
Further, we prove that the condition iv) ( Definition 4.3) holds
if and only if  $\Phi$ is maximal isotropic subspace and 
$\Phi \oplus\t=\h$.\\
2) Our goal in 2) is to prove that $\Phi\bigcap\t\ne 0$ is equivalent
to $\phi$ has no fixed points.

Let $\Phi\bigcap\t\ne 0$ then there exists $x\in\h$ such that $\qq x=x$
and $x\in\Phi\subset\C\Phi=\lh_1$. 
Hence, $x ={\qq x}\in{\qq\Phi}\subset\C{\qq \Phi}=\lh_2$.
It follows that $x\in\lh_1\bigcap\lh_2$.

Denote by $x_1$ the image of $x$ in $\lh_1/\l_1$ and by     
$x_2$ the image of $x$ in $\lh_2/\l_2$. Note that, if $x\in\Phi$, then
$^\ll x_1=x_1$, and ,if $x\in\t$ , then $\qq x = x$ , $\qq x_1=x_2$.
We obtain $\phi(x_1) = {^{\q\ll}x_1} = {\qq x_1} =x_2$. Therefore,
$\phi(x_1)=x_2$ and $\phi$ has fixed point.

On the other hand, let $\phi$ has fixed point ( i.e. there exists $x\ne 0$,
$x\in\lh_1\bigcap\lh_2$ such that $\phi(x_1)=x_2$).
 Consider complex subspace 
$$\h_\phi = \{x\in\lh_1\bigcap\lh_2\vert \phi(x_1)=x_2\}$$
We are going to show that $\h_\phi$ is $\q$-invariant. Indeed, if
$x\in\h_\phi$, then $\phi(x_1)=x_2$,  ${^\q\ll}x_1=x_2$
and ${^\ll x_1}={\qq x_2}$. Note, that $({^\q x})_1 = {^\q x} +\l_1=
{^\q (x+\l_2)}= {^\q x_2}$. Thus, ${^\ll x_1} = ({^\q x})_1$ and
$$\phi({^\q x})_1 = {^{\q\ll}({^\q x})_1} = {^{\q\ll\ll}x_1} = ({^\q x})_2$$
This proves that $\h_\phi$ is $\q$-invariant.

There exists $x\ne 0$, $x\in\h_\phi$, ${^\q x}=x$. The last means that
$x\in\t$. We shall show that $x\in\Phi$.
The $\phi(x_1)=x_2$ implies ${^{\q\ll}x_1}=x_2$, ${^\ll x_1}={^\q x_2}=
({^\q x})_1= x_1$. Finally, ${^\ll x_1}= x_1$, $x_1\in\Phi/\l_1$ and
$x\in\Phi$. Hence, $\Phi\bigcap\t\ne 0$. \\
3) We shall show that $\phi $ is isotropic with respect to $Q_\R$ if and only
if $\phi$ save the $K(\cdot,\cdot)$.

Let $\Phi$ be an isotropic subspace with respect to $Q_\R$. 
This means that $K(x,y)= K({^\q x}, {^\q y})$ for all $x,y\in\Phi$.
Recall that $\l_1$ is the greatest complex subspace in $\Phi$.
For $x\in\l_1$ , $y\in\Phi$ we have

$$\left \{ \begin{array}{rcl}K(x,y)&=& K({^\q x}, {^\q y})\\
K(ix,y)&=& K({^\q ix}, {^\q y})\end{array}\right.$$
It follows $iK(x,y)= -iK(x,y)$ and $K(x,y)= 0$.
We prove that $\l_1\subset Ker(K\vert_\Phi)$.

Notice that for $x\in\Phi$ the following equalities hold
$ {^\ll x_1}=x_1$, ${^\q \phi (x_1)}= x_1$, $\phi(x_1)={^\q x_1}= ({^\q x})_2$.
 Then for all $x,y\in\Phi$ we obtain
$ K(x_1,y_1)= K(x_1+\l_1,y_1+\l_1) = K(x,y) = K({^\q x}, {^\q y}) =
 K({^\q x}+\l_2, {^\q y}+\l_2) = K(({^\q x})_2, ({^\q y})_2) = 
K(\phi(x_1), \phi(y_1))$.
Therefore, $\phi:\lh_1/\l_1\mapsto\lh_2/\l_2$ save the Killing form.

On the other hand, if $\phi$ save $K(\cdot, \cdot)$, then
$\l_1\subset Ker(K\vert_\Phi)$ and $\l_2\subset Ker(K\vert_{\qq\Phi})$.
We have
$$ K(x_1,y_1)= K(x_1+\l_1,y_1+\l_1) = K(x,y)$$
$$ K(\phi(x_1), \phi(y_1)= K(({^\q x})_2, ({^\q y})_2) 
= K({^\q x}+\l_2, {^\q y}+\l_2)= K({^\q x}, {^\q y})$$
We see $K(x_1,y_1)= K({^\q x}, {^\q y})$ for all $x,y\in\Phi$
and $\Phi$ is isotropic subspace. This proves 3).

Thus , $\Phi$ is isotropic subspace and $\Phi\oplus\t=0$ if and only if
$\phi$ satisfies the conditions i)-v) with $\1=\2-\emptyset$.\\
4) Using $\phi$, we can construct a complex Manin triple with
$W= \W$ ( see Theorem 1.4). We are going to show that 
$\W^\Gamma=\Phi\oplus\n^+$.

Every $x\in\W$ admites the presentation  
$$ x = l_1 e + n_1^+ e + he + \phi(h)f+ n_2^ - f + l_2 f\eqno(4.29),$$ where
$ l_i\in\l_i$, $n_1^+\in\n_1^+$, $h\in\lh_1$, $n_2 \in\n_2^-$ .
Then $ {\qq }x = {\qq l_1} f + {\qq n_1^+}f + {\qq h}f + 
{\qq \phi(h)}e+ {\qq n_2^ -}e + {\qq l_2}e$.

Recall that ${\qq l_1}= \l_2$ and ${\qq l_2}= \l_1$.
The element $x$ is $\q$-invariant (i.e.${\qq x}=x$).  
if and only if ${\qq n_1^+}= n_2^-$ ,
${\qq l_1} = l_2$
and ${\qq \phi(h)}=h (mod(\l_1))$.
The last means that ${^\ll h}= h(mod(\l_1))$ and $h\in\Phi$.

The $\W^\Gamma$ is spanned by elements 
$n_1^+e +{\qq n_1^+}f$, $l_1e +{\qq l_1}f$,
$h e+ {\qq h}f$, $h\in\Phi$. 
As we see above the double algebra $\dR = (\dC)^{\Gamma} $ 
consists of the
elements $\tilde x = xe+{^\q x}f.$ 
The map $\tilde x\mapsto x$  is an $\R$-isomorphism of $\dR$ to  $\gC$.
The image of $\W^\Gamma$ coinsides with $\Phi\oplus\n^+$.\\
5) We are going to prove that 
$(Ad_g\W)^\Gamma= Ad_g(\Phi\oplus\n^+)$, $g=r\o$, $r\in G(\R)$, $\o\in\O$.

Recall ${^\q\o}=\o a$, $a\in H$ and $a^2=1$.
The algebra $Ad_g\W$ equals to
$$Ad_g(\l_1)e + Ad_g(\n^+) e +
\sum_{x\in\lh/\l_1}\C(Ad_g(x)e + Ad_g(\phi(x))f)
+ Ad_g(\n^-)f+ Ad_g(\l_2)f$$
We have 
${^\q Ad_g\n^+}$ = ${^\q Ad_r Ad_\o\n^+}$ = $Ad_rAd_{{^\q\o}}(^\q\n^+)$=
$Ad_rAd_{{^\q\o}}Ad_a(\n^-)$ =
$Ad_r Ad_\o(\n^-)$ = $Ad_g \n^-$
and ${^\q Ad_gh}$ = ${^\q Ad_r Ad_\o h}$ = $Ad_rAd_{{^\q\o}}(^\q h)$ =
$Ad_rAd_{{^\q\o}}Ad_a(h)$=
$Ad_rAd_\o(h)$ = $Ad_g h$ for all $h\in\h$.
It follows that ${^\q (Ad_g(\l_1)} = Ad_g(\l_2)$,
${^\q (Ad_g(\lh_1)} = Ad_g(\lh_2)$.

Let ${\qq(Ad_gx)}=Ad_gx$, $x\in\W$ and presented as in (4.29). Then
${^\q (Ad_gn_1^+)} = Ad_g(n_2^-)$, ${^\q (Ad_g(l_1)} = Ad_g(l_2)$,
${^\q Ad_g\phi(h)} = Ad_gh(mod(\l_1))$. Since $\phi(h)\in\h$,
 then ${^\q Ad_g\phi(h)} = Ad_g{^\q\phi(h)} = Ad_g{^\ll h}$. 
We get $Ad_g{^\ll h}= Ad_gh$. Thus, ${^\ll h = h} $ and $h\in\Phi$.

The algebra of invariants $(Ad_g\W)^\Gamma$ is generated by 
$Ad_gn_1^++ {\qq (Ad_gn_1^+)}f$, $Ad_gl_1 +{\qq (Ad_gl_1)}f$,
$Ad_gh + {\qq (Ad_gh)}f$, $h\in\h$.
Hence, $(Ad_g\W)^\Gamma= Ad_g(\Phi\oplus\n^+)$.
This completes the proof. $\Box$\\
{\bf Corollary 4.6} Every class of weak equivalent Manin triples 
$(\gR,W(\R),\gC)$ with the inner form $\gR$ consists
of finite number of gauge equivalent Manin triples.$Box$\\
The proof followes form the finiteness of $\O$.

Our next goal is to get the classification of Manin triples 
for outer forms. \\
Notations:\\
$\0=\1\bigcup\2$;\\
$\hP=\{v\in\h\vert \a(v) = \phi(\a){v}, \a\in\1\} = 
\{v\in\h\vert [v,\phi(x)] = \phi([v,x]), x\in\lr_1\}$;\\
The last formula can be rewritten as $ad_v\phi(x)=\phi ad_v(x), x\in\lr_1$
We shall call $\hP$ the stabilizer subalgebra of $\phi$.\\
{\bf Lemma 4.7}. Let $s$ be as in (3.1), $s(\1)=\2$ and 
$^s\phi(\a)=\phi^{-1}s(\a)$ ( in other words $({^s\phi})^2(\a)= \a$)
for all $\a\in\1$.
Then subalgebra $\hP$ is $\q$- invariant.\\
{\bf Proof}. First we shall note  that if $u=\sum_{i=1}^n\mu_iH_i\in\h$ then
$\a({^\q(u)})= -\sum_{i=1}^n\overline{\mu_i}\a(s(H_{i}))$
We have 
$$\a({^\q u})=-\overline{^s\a(u)}$$
Applying this equality, we see:
$$\phi(\a)({^\q v})=-\overline{^s\phi(\a)(v)}$$
Let $v\in\h^\phi$. Since, $s(\phi(\a))\in\1$, then
$$-\overline{^s\phi(\a)(v)} = -\overline{\phi{^s\phi(\a)(v)}} =
-\overline{\phi(\phi^{-1}(s(\a)))(v)}=
-\overline{^s\a(v)}=\a({^\q v})$$ 
 This proves Lemma 4.7.$\Box$\\
{\bf Corollary 4.8}.The subalgebra $\hP$ admits a real form
$\hPR =\{x\in\h\vert {^\q x}=x\}$.\\
{\bf Definition 4.9}([BD]). Let $\phi:\1\mapsto\2$ satisfy i),ii),iii). 
Let $\a,\b\in\1$. We shall say that $\a < \b$ if there exists integer
$k$ such that $\phi^k(\a)=\b$.\\
Note that $\a < \phi(\a)<\cdots<\phi^k(\a)=\b$.\\
{\bf Definition 4.10}. The set $\a_1<\a_2<\cdots<\a_k$ is called a chain.
We shall denote by $C(\a)$ the maximal chain with $\a$.\\
{\bf Lemma 4.11}. 1) The maximal chains either coinside or don't intersect.;
2) If $\phi$ and $s$ as Lemma 4.7 then the image $s(C)$ of a maximal chain
is also a maximal chain.\\
{\bf Proof}.\\
1) Let $C_1= \{\a_1 < \a_2 <\cdots<\a_k\}$ and 
$C_2= \{\b_1 < \b_2 <\cdots<\b_k\}$ are maximal chains. If
$\phi^n(\a) = \phi^m(\b)$ , $n\ge m$ then $\phi^{n-m}(\a)=\b$ and
, therefore, $\b\in C_1$. Hence, $C_2\subset C_1$. Since $C_2$ is maximal,
 $C_1 =C_2$.\\
2) Let $\{\a_1 < \a_2 <\cdots < \a_k\}$ be a maximal chain with 
$\a_{i+1} = \phi(\a_i)$. Denote $\b_i=\phi(\a_i)$. We see
$\phi(\b_{i+1}) = \phi(s\a_{i+1})) = \phi(s(\phi(\a_i))) = s(\a_i) =\b_i$.
Hence, $\b_1 < \b_2 <\cdots < \b_k$.$\Box$\\
{\bf Corollary 4.12} $\0$ is decomposed in maximal chains. The automorphism $s$
acts on the set of maximal chains.\\ 
Notations:\\
$\HP=\{p\in H\vert Ad_p\phi(x) = \phi Ad_p(x)\}$. We shall call $\HPR$ the 
stabilizer of $\phi$. The Lie algebra of $\HP$ is $\hP$.\\
$\HPR$ is Lie subgroup of the subalgebra $ \hPR$.;\\
For $p=exp(v)$, $v\in\h$ denote $p_a=exp(a(v))$ (see Proof of Theorem 4.4);\\
Denote by $\APR$ the subgroup of $\HPR$ consisting of the elements
$a$ such that\\
 1) $a^2=1$ i.e. $a_\a=\pm1$ for all $\a\in\Pi$,\\
2) $a_\a=1$ for 
$\a\in\Pi-\0$, $s(\a)\ne\a$ and $\a\in \0$, $sC(\a)\ne C(\a)$.;\\
$\APR' = \{a\in\APR\vert\quad\exists g\in G, {\qq g^{-1}}g=a\};$\\
$\O(\APR) = \{g\in G\vert\quad{\qq g^{-1}}g \in \APR\}$;\\
The group $G(\R)$ acts on $\O(\APR)$  by $g\mapsto rg$, $r\in G(\R)$.
Denote by $\OP$ the set of representatives of $G(\R)$- orbites in $\O(\APR)$.
The set $\OP$ is finite.\\
{\bf Theorem 4.13} Let $\gR$ be an outer form of $\g=\gC$ . The conjugation
of $\gR$ as in (3.1) with $\e=1$. 
The map $\phi $ as in Theorem 1.4. Then\\
1) $\W $ is $\q$-invariant if and only if and
$\phi$ satisfies v).\\
2) There exists a $\q$-invariant subalgebra in $Ad_{G(\C)}\W$ iff
$\phi$ satisfies v). Every $\q$-invariant subalgebra in 
$Ad_{G(\C)}\W$ has a form $Ad_r Ad_\o\W$, where $r\in G(\R)$, $\o\in\OP$.\\
{\bf Proof}. 
The proof is divided into several steps.
We recall that (4.1) implies (4.2) with $p={^\q g^{-1}}g$. 
In the first step we shall show that, if $\W$ is $\q$-invariant then
the condition v) holds. In the steps 2)-3) we will see that the above element 
$p$ commutes with $\phi$ and belongs to Cartan subgroup.
The statements of the Theorem will be proved in steps 4)-5).\\
1) 
We are going to show in this step that
if $Ad_g\W$ is $\q$-invariant then $\phi$ satisfies v).

As we saw in (4.26) $^\q\lr_1=\lr_2$. 
It was proved that (4.1) implies (4.4). We have 
$(^\q Ad_{p_2}\phi)^2(x)= x$ for $x\in\lr_1/\l_1$.\\
Decompose $\1$ and $\2$ in connected components:
$$\1 = \bigcup_{i=1}^k \Pi_{1i}$$
$$\2 = \bigcup_{i=1}^k \Pi_{2i}$$
Then $\g_1$ and $\g_2$ are decomposed in direct sums:
$$\g_1 = \sum_{i=1}^k \g_{1i}$$
$$\g_2 = \sum_{i=1}^k \g_{2i}$$
We may assume that $\phi(\g_{1i})=\g_{2i}$.\\
{\bf Notaton}. If $\mu:\g_1\mapsto\g_1$ is a automorphism 
(linear or semilinear) then we shall denote by $\hat\mu$ the 
corresponding substitution of the set $\{ 1,2,\cdots,k\}$ of 
connected components.\\
One can see that 
${^\q Ad_{p_2}\phi}$ and $ {^\q\phi}$ yield the same 
substitution $\widehat{\qq \phi}$
of components of the subalgebra $\g_1$. 
Hence the substitutions of maps 
$({^\q Ad_{p_2}\phi})^2 = 1$ and  $(\qq\phi)^2$ are also equal.
then
$$(\widehat{{^\q\phi}})^2 = 1$$
Thus, $ ({^\q\phi})^2$ is $\C$-linear automorphism of $\g_1$, 
saving components.
Denote by $\phi_i$ the restrition of $\phi$ on $\Pi_{1i}$.
Save this notation for the restriction of $\phi$ on $\g_{1i}$. 
Since $(\widehat{\qq\phi_i})^2 = 1$, then there two possibilities 
$({\qq\phi_i})^2 =1$ or $({\qq\phi_i})^2$ equals to some 
nontrivial automorphism of $\Pi_{1i}$ . We shall denote it 
by $ s_{1i}$.

 Suppose that 
$({^s\phi_i})^2=s_{1i}$ for some $i$.
Hence,
$${^\q\phi_i}{^\q\phi_i}=s_{1i}$$
Denote $\psi_i={^\q\phi_i}^\q : \g_{2i}\mapsto\g_{1i}$.
We get $\psi_i\phi_i=s_{1i}$ is an outer automorphism of $\\g_{1i}$.\\
On the other hand, (4.4) implies 
$${^\q Ad_{p_2} \phi_i}{^\q Ad_{p_2}\phi_i}=1$$
Then
$$Ad_{p_1^{-1}} \psi_i Ad_{p_2}\phi_i=  
Ad_{p_1^{-1}}{^\q\phi_i}^\q Ad_{p_2}\phi_i = 
({^\q Ad_{p_2}\phi_i})^2 = 1$$
The above equality we rewrite in the form
$$Ad_{p_1^{-1}}\psi_i\phi_i(\phi_i^{-1} Ad_{p_2}\phi_i) = 1\eqno (4.30)$$
We are going to show that $\phi^{-1} Ad_{p_2}\phi_i$ is inner automorphism of
$\g_{1i}$.

Really, let $p_2$ be decomposed $p_2= g_2h_2$, $g_2\in G_2$, $h_2\in H$ and
$g_2h_2= h_2g_2$. There exists $v\in\g_2$ such that $g_2= exp(v)$ . 
Denote $u = \phi^{-1}(v)$.
The equality $\phi^{-1}[v,\phi(x)]=[\phi^{-1}(v),x] = ad_u(x)$, $x\in\g_1$
 yields $\phi^{-1}ad_v\phi = ad_u$ and 
$\phi^{-1}Ad_{p_2}\phi = Ad_q$, $q=exp(u)$.
Substitute this to (4.30): $Ad_{p^{-1}}\psi_i\phi_i Ad_q = 1$.
Therefore, $\psi_i\phi_i$ is inner automorphism of $\g_{1i}$. This leads to 
contradiction. 
Finally, $({\qq\phi_i})^2 =1$ for all $i$. This proves v). The step 1) is 
completed. \\
2) As we see in 1)  the (4.4) and (4.28) (condition v) )
holds samultaneously. Therefore, we have
$$ x= (\qq Ad_{p_2}\phi)^2(x)= \qq Ad_{p_2}\phi\qq Ad_{p_2}\phi(x)=
 Ad_{p_1^{-1}}{\qq\phi\qq} Ad_{p_2}\phi(x)= 
Ad_{p_1^{-1}} \phi^{-1}Ad_{p_2}\phi(x)$$ 
for all $x\in\lr_1/\l_1$.
It follows
 $$Ad_{p_2}\phi(x) = \phi Ad_{p_1}(x),\eqno(4.31)$$ for $x\in\lr_1/\l_1$.
This proves 2).\\
3) In this step we are going to show that $p_1=p_2\in H^\phi$. \\
There exists $v_i\in\r_i$ such that $p_i= exp(v_i)$.
According Lemma 4.2, we have
$v_1=\mu_1 +\th + t$, 
$v_2= \mu_2+\th+t$ where $\mu_i\in\m_i$, 
$\th\in(\g_1\bigcap\g_2)\bigcap(\n^+\oplus \n^-)$, $t\in \h$.
We shall prove further that $\mu_1=\mu_2=\th=0$.
 
We can decompose $t=t_i+t_i'$ where $t_i\in\h_i=\g_i\bigcap\h$, $t_i'\in\h$
and $[t_i',\g_i]= 0$. Denoting $g_i=exp(\mu+\th+t_i)$, $h_i=exp(t_i')$, we get
$p_i= g_ih_i= h_ig_i$. 

The (4.31) implies 
$\phi(Ad_{g_1}x) = Ad_{g_2}\phi(x)$. 
Let $G_i$ be the subgroup of $G$ with Lie algebra $\g_i$.
We can extend $\phi$ to the isomorphism of
the groups $\phi:G_1\mapsto G_2$. Then 
$Ad_{\phi(g_1)}\phi(x) = Ad_{g_2}\phi(x)$. Hence, $\phi(g_1)=g_2$ .
Therefore, $\phi(\mu_1+\th + t_1)= 
\mu_2+ \th + t_2$. The element $\th$ is a sum $\th= \th^++ \th^-$ where
$\th^\pm\in\g_1\bigcap\g_2\bigcap\n^\pm$. We obtain
$$\phi(\mu_1)+\phi(\th^+) + \phi(\th^-) + \phi(t_1)= 
\mu_2+ \th^++ \th^-+ t_2\eqno(4.32)$$
Recall that the map $\phi$ is defined via the map $\phi:\1\mapsto\2$.
Thus, $\phi(\n^+\bigcap\g_1)\subset\n^+\bigcap\g_2$, 
$\phi(\n^-\bigcap\g_1)\subset\n^-\bigcap\g_2$, $\phi(\h_1)\subset\h_2$.
The equality (4.32) holds if and only if
$$\left \{ \begin{array}{rcl}\phi(\th^+)&=& \mu_2+\th^+\\
 \phi(\mu_1)+\phi(\th^-)& = &\th^-\\
\phi(t_1)&=&t_2\end{array}\right.\eqno (4.33)$$
Suppose that $\th^+\ne0$. Denote $\Delta_i^+$ the system if positive roots
generated by $\Pi_i$. Then
$$
\th^+=\sum_{\a\in\Delta_1^+}\xi_\a E_\a,\qquad
\mu_2= \sum_{\a\in\Delta_2^+-\Delta_1^+}\xi_\b E_\b$$
The equality $\phi(\th^+)=\mu_2+ \th^+$ implies
$$
\sum_{\a\in\Delta_1^+}\xi_\a E_{\phi(\a)} = 
\sum_{\a\in\Delta_2^+-\Delta_1^+}\xi_\b E_\b
+ \sum_{\a\in\Delta_1^+}\xi_\a E_\a\eqno(4.34)
$$
For every simple root $\a\in\1$ we denote $l(\a)$ the length 
of the maximal chain  $C(\a)$ (see definition 4.10) from the 
beginning to $\a$. 
Denote
$$l_0=min\{l(\a)\vert \a\in\Delta_1^+, \xi_\a\ne0\}$$
and $l_0=l(\a_0)$ for some positive root $\a_0\in\Delta_1^+$
The term $\xi_{\a_0}E_{\a_0}$ belongs to the right side of (4.34)
and dose not belong to to left side. A contradiction.

Therefore $\th^+=0$ and $\mu_2=0$.
Similarly we get  $\th^-=0$ and $\mu_1=0$.
This proves $v_1=v_2=t$ , $p_1=p_2=exp(t)=p\in H$
and 
$$ Ad_p\phi = \phi Ad_p\eqno(4.35)$$
4) So $p=exp(v)\in\HP$, $v\in\hP$.
As in the proof of Theorem 4.4., we shall identify $p$ with the
system $(p_{\a})_{\a\in\Pi}$, $p_\a=exp(\a(v))$.
Since $v\in\hP$, then $p_{\phi(\a)} = p_\a$ and
$p_\a$ is constant on the maximal chains in $\0$.

The condition ${^\q p}=p^{-1}$ is equivalent to 
$\overline{p_{s(\a)}}=p_\a$. If $s(\a)=\a$ then it means 
$\overline{p_\a}=p_\a\in\R$.\\
In the case $p={^\q h}^{-1}h$ the system will be
$p_\a=\overline{h_{s(\a)}}h_\a$.
If $s(\a)=\a$, then $p_\a=\overline{h_\a}h_{\a} > 0$.
If $s(\a)=\a'\ne\a$, then $p_\a =\overline{h_{\a'}}h_\a$ and
$p_{\a'} =\overline{h_{\a}}h_{\a'}$.\\
Using Lemma 4.11, one can prove that for every $p\in\HP$ there exists
$\h\in\HP$ such that  the element $a=h^{\phi} ph^{-1}$ satisfies $a^2=1$.
The  element $a\in\HP$ and we can choose $h$ so that $a_\a=1$ for all
$\a\in\Pi-\0$, $s(\a)\ne\a$ and $\a\in \0$, $sC(\a)\ne C(\a)$.\\
Since ${^\q a} = a^{-1}$ and $a^2=1$, then ${^\q a} = a$ and so $a\in\HPR$.\\

We denote by $\o$ in $\OP$ such that
$ a = {^\q}hph^{-1} ={^\q}h{{^\q}g^{-1}}gh^{-1}= 
{^\q \o^{-1} \o}.$ 
  $$ ^\q (gh^{-1}\o^{-1})^{-1} gh^{-1}\o^{-1} = 1$$
Therefore, $gh^{-1}\o^{-1} = r \in G(\R)$
and $g=r\o h$.
We have
$$Ad_g\W = Ad_rAd_{\o}Ad_h\W =  Ad_r Ad_{\o}\W $$
Thus, we have proved that , if subalgebra $Ad_g\W$ is $\q$-invariant , then
$\phi$ satisfies i)-iv)) and $g=r\o$, 
$\o\in\O$. Putting $g=1$, we get 1).\\
5) Now we shall prove the converse statement : all this subalgebras are 
$\q$-invariant.
Recall
$$Ad_g\W = Ad_g\l_1e+Ad_g\n_1^+e+
\sum_{x\in\r_1}\C\left(Ad_gxe+Ad_g\phi(x)f\right) + 
Ad_g\n_2^-f+ Ad_g\l_2f$$
If $g=r\o$, $r\in\GR$, $\o\in\OP$ i.e ${^\q\o}= \o a$ and $a\in \APR'$
We have
${^\q (Ad_g\n_1^+)}$ = ${^\q (Ad_rAd_\o\n_1^+)}$=
$Ad_rAd_{^\q\o}{^\q\n_1^+}$ = $Ad_rAd_\o Ad_a\n_2^-$ = $Ad_g\n_2^-$ .
Similarly, 
${^\q (Ad_g\n_2^+)}$ = $Ad_g\n_1^+$
and ${^\q (Ad_g\l_1)}$ = ${^\q (Ad_rAd_\o\l_1)}$=
$Ad_rAd_{^\q\o}{^\q\l_1}$ = $Ad_rAd_\o Ad_a\l_2$ = $Ad_g\l_2$ .

Recall that $a\in\HPR$ and $\HP$ is stabilizer of $\phi$.
We have
$${^\q Ad_g\left(xe+\phi(x)f\right)} = 
{^\q Ad_rAd_\o\left(xe+\phi(x)f\right)} =
Ad_rAd_{^\q\o}\left({^\q x}f+{^\q\phi(x)}e\right) =$$
$$Ad_rAd_\o\left(Ad_a{^\q x}f+Ad_a{^\q\phi(x)}e\right) =
Ad_g\left({^\q Ad_ax}f+{^\q\phi(Ad_ax)}e\right) =$$
$$Ad_g\left(\phi^{-1}{^\q(Ad_ax)}e+{^\q Ad_ax}f \right).$$
Since $\phi^{-1}{^\q(Ad_ax)}\in\r_1$ 
and ${\qq Ad_gx}= \phi(\phi^{-1}{^\q(Ad_ax)})$,
then proves that $\W$ is $\q$-invariant.
This completes the proof of Theorem 4.3.$\Box$\\
{\bf Theorem 4.14}. Let $\gR$ be an outer form of $\gC$ with conjugation
(3.1) with $\e=1$.
Let $\phi:\1\mapsto\2 $ be a map satisfying i)-v).
For every $a\in\APR'$ consider the conjugation 
$\ll_a(x)=Ad_a{^\q\phi(x)}$ on $\lr/\l_1$.
By $\FPa$ we shall denote  the real form 
$\{x\in\lr_1/\l_1\vert {\ll_a(x)} = x\}$
in $\lr_1/\l_1$.In the case $a=1$ put $\FP=\FPa$. Denote 
$$\WPa = \FPa+\l_1+\n_1^+$$
Then:\\
1) $(\gR,\WPa,\gC)$ is a Manin triple for all $a\in\APR'$. 
In particular, for if $a=1$ then $(\gR,\WP, \gC)$ with $\WP=\FP+\l_1+\n_1^+$ 
is a Manin triple;\\
2) Every Manin triple $(\gR,W,\gC)$ is weak equivalent to
 some triple $$(\gR, \WPa, \gC);$$
3) Every Manin triple $(\gR,W,\gC)$ is gauge equivalent to the Manin triple
$$(\gR,Ad_\o\WPa,\gC),$$ where ${^\q\o}^{-1}\o=a$ and $\o\in\OP$.\\
{\bf Proof.} As in Theorem 4.4 
the double $\dR$ coinsides with $(\dC)^{\G}$ and 
$\R$-spanned by $\tilde x = xe+{^\q(x)}f$, $x\in\gC$.
We are going to study algebra of invariants $(Ad_g\W)^\Gamma$.
as we have got in Theorem 4.13 ${\qq Ad_g\n_1^+}= Ad_g\n_2^-$,
${\qq Ad_g\l_1}= Ad_g\l_2$ and the algebra 
$\{Ad_g(xe+\phi(x)f)\}$, $x\in\lr_1/\l_1$  is $\q$-invariant.
 
The $\q$-invariant subalgebra $(Ad_g\W)^\G$ is $\R$-spanned by the elements
$Ad_gn_1^+e+{^\q Ad_gn_1^+}f$, $n_1^+\in\n_1^+$ ;
$Ad_gl_1e+{^\q Ad_gl_1}f$, $l_1\in\l_1$ 
and
$$ Ad_gxe + Ad_g\phi(x)f + {^\q(Ad_gxe + Ad_g\phi(x)f)} =$$
$$(Ad_gx+ {^\q Ad_g\phi(x)})e + ({^\q Ad_gx} + Ad_g\phi(x))f=$$
$$(Ad_gx+ {^\q Ad_g\phi(x)})e + {^\q(Ad_gx+ {^\q Ad_g\phi(x)})}f,$$
with $x\in\lr_1/\l_1.$\\
The map $\tilde x\mapsto x$ establishes 1-1 correspodence between $\dR$ and
$\gC$. The image of $\q$-invariant subalgebra $Ad_g\W$ in $\gC$ is
$\R$-spanned by the elements
$Ad_g\n_1^+$, $Ad_g\l_1$ and 
$Ad_gx+{\qq Ad_g\phi(x)}$= 
$Ad_rAd_\o(x)+{\qq Ad_rAd_\o \phi(x)}$=
$Ad_rAd_\o(x)+Ad_rAd_{\qq\o}{\qq\phi(x)}$=   
$Ad_rAd_\o(x)+Ad_rAd_\o Ad_a{\qq\phi(x)}$   
= $Ad_rAd_\o(x+\ll_a(x)),$
$x\in\lr_1/\l_1$.
So the image coinsides with $\WPa = Ad_rAd_\o(\FPa+\l_1 + \n_1^+)$. \\
{\bf Corollary 4.15} Every class of weak equivalent Manin triples 
$(\gR,W(\R),\gC)$ consists
of finite number of gauge equivalent Manin triples.\\
The proof followes form finiteness of $\O$ for inner forms
and from finiteness of $\OP$ for outer forms.\\
\begin{center}
{\bf Examples}
\end{center}
We shall illustrate the Theorems 4.5 and 4.14 by the examples.\\
{\bf Example 1}.$\d = \gC = sl(2,\C)$, $\gR = su(2)$\\
Let $H, E^+, E^-$  be a Cartan basis of $sl(2,\C)$.\\ 
All Manin triples with $\g=su(2)$ and $\d = sl(2,\C)$
have the form $(\g,V,\d)$ with 
$$V = Ad_r(\R H + \C E^+) = Ad_r\left(\begin{array}{cc}a&z\\0&-a\end{array}
\right ),$$ 
where $a\in\R$, $z\in\C$ and 
$$r\in G(\R)= \left(SL(2,\C)/\pm E\right)^\Gamma= Ad_{SU(2)}.$$
{\bf Example 2}. Consider the algebra 
$$\gR= \left(\begin{array}{cc}\a&\b\\ \bar\b&\bar\a\end{array}
\right )$$ 
with $\a, \b\in\C$, $\a+\bar\a=0$. This algebra is isomorphic to $sl(2,\R)$.
We get $\gC=sl(2,\C)$ ,
$$G(\C)= SL(2,\C)/\pm E\quad,\quad G(\R)= 
\left(SL(2,\C)/\pm E\right)^\Gamma.$$ 
The Cartan subgroup $H(\R)=diag(\a,\bar\a)$ with $\vert\a\vert=1$.
The subgroup $\AR$ of $G(\R)$ consists of two elements
$\pm E, \pm diag(i,-i)$.
By definition, 
$\AR'$ consists of $a\in\AR$ such that there exists 
$g\in G, {\qq g^{-1}}g=a$. Simple calculatons
yield $\AR'=\pm E$ and $\O=\pm E$.
All Manin triples $(\gR,V,\gC)$ are gauge equivalent to 
$\R H + \C E^+.$ 

The $G(\R)$ is presented by $g\in SL(2,\C), {^\q g}=\pm g$.
This elements lies in the union $SL(2,\R) \bigcup SL(2,\R)w_0$ where
$$w_0= \left (\begin{array}{cc}0&1\\-1&0\end{array}
\right ).$$
We get more precise answer : all Manin tripes $(\gR,V,\gC)$ are 
$SL(2,\R)$- equivalent to 
$\R H + \C E^+$ or to 
$\R H + \C E^-$. One can prove that this Manin triples are not 
$SL(2,\R)$-equivalent.\\
{\bf Example 3}. In this example we treat the $\gR = sl(2,\R)$ in the usial
presentation.
Let $$u_0=\frac {1}{\sqrt 2}\left (\begin{array}{cc}1&i\\i&1\end{array}
\right )$$
Every Manin triple with $\g=sl(2,\R)$, $\d=sl(2,\C)$ 
has a form $(\g,V,\d)$ with  $V=Ad_rAd_u (R H+\C E^+)$ where  $r\in SL(2,\R)$,
$u={u_0}$ or $ u= \overline{u_0}$. There exist exactly two 
bialgebra structures on $sl(2,\R)$ up to $SL(2,\R)$-equivalence.\\
{\bf Example 4}. $\d = \gC = sl(3,\C)$,
$$\gR = \left (\begin{array}{ccc}x&y&z\\v&t&{\overline v}\\
{\overline z}&{\overline y}&{\overline x}\end{array}\right ) $$
 with the entries $x,y,z,v\in \C$, $t\in \R$ and $x+{\overline x}+t=0$.
 This algebra is isomorphic to $sl(3,\R)$.\\
1) $\Pi = \{\a_1,\a_2\}$, $\phi:\a_1\mapsto\a_2$,
$\r_1=\g_1+\h$, $\r_2=\g_2+\h$.\\
There exists a unique extention of $\phi$ to $\phi: \r_1\mapsto \r_2$
satisfying iv) i.e. $\phi\vert_\h$ preserve the Killing form and has no
fixed points  $\phi( diag(x_1,x_2,x_3)) = diag(x_3,x_1,x_2)$.
By calculations, $$\APR =\{diag(1,1,1), ~diag(-1,1,-1)\}= \APR'$$
$$\OP= \{diag(1,1,1),~diag(-1,-1,1)\}.$$
Using Theorem 4.14, we see that there exist two types of 
Manin triples, corresponding to the map
$\phi:\a_1\mapsto\a_2$.
In every  Manin triple $(\gR,W,\gC)$ the complementary subalgebra $W$ 
has one of the forms
$$ Ad_r\left (\begin{array}
{ccc}\a&\b&z_1\\{\overline \b}&{\overline{\a}}&z_2\\
0&0&-2Real(\a)\end{array}\right );
\qquad\qquad
Ad_r\left (\begin{array}
{ccc}\a&\b&z_1\\-{\overline \b}&{\overline{\a}}&z_2\\
0&0&-2Real(\a)\end{array}\right )$$
with the entries $\a, \b, z_1, z_2\in \C$ and 
$$r\in G(\R)=\left(SL(3,\C)/Z\right)^\Gamma.$$ The center 
$Z= \{1, \epsilon, \epsilon^2\}$, $\epsilon^3= 1$. 
Here  $G(\R)$ is presented by $g\in SL(3,\C), {^\q g}=\epsilon^i g$.
This elements lies in the union $SL(3,\R) \bigcup SL(3,\R)\epsilon
\bigcup SL(3,\R)\epsilon^2 $. We see that $G(\R)=Ad_{SL(3,\R)}.$\\
2) The case $\Pi = \{\a_1,\a_2\}$, $\phi:\a_2\mapsto\a_1$ is treated 
similary  to 1). We get the following $W$:\\
$$ 
Ad_r\left (\begin{array}
{ccc}-2Real(\a)&z_1&z_2\\ 0&\a&\b\\0&{\overline \b}&{\overline{\a}}\\
\end{array}\right );
\qquad\qquad 
Ad_r\left (\begin{array}
{ccc}-2Real(\a)&z_1&z_2\\ 0&\a&\b\\0&-{\overline \b}&{\overline{\a}}\\
\end{array}\right )
$$
where $r\in SL(3,\R)$.\\
3) $\Pi_1=\Pi_2=\emptyset$.\\
Let $\t = \h\bigcap\gR$ and
$\Phi$ is a real subalgebra of $\h$ such that $\Phi$ is isotropic 
with respect to $Q_\R(x,y)=2iK(x,y)$ and $\Phi\oplus \t=\h$.
Every Manin triple $\gR,W,\dR$ is $Ad_{SL(3,\R)}$ equivalent to Manin 
triple with $W= \Phi+\n$.

We get the classification of Manin triples over $sl(3,\R)$
up to $SL(3,\R)$-equivalence.

\end{document}